\documentclass[12pt]{article}

\newtheorem{theorem}{Theorem}[section]
\newtheorem{corollary}{Corollary}[section]

\def\square{\hbox{\vrule\vbox{\hrule\phantom{o}\hrule}\vrule}}
\def\endofproofsymbol{\square}

\newenvironment{proof}{\par\noindent\textbf{Proof.}\hskip\labelsep}
                      {\unskip\ \nobreak\hbox{}\hfill
                       \endofproofsymbol{\parfillskip=0pt\par}
                       \vspace{\baselineskip}}

\def\eqref#1{\mbox{(\ref{eq:#1})}}    
\def\BBR{{\it I\kern-0.408em R}}
\makeatletter \def\@cite#1#2{{#1\if@tempswa , #2\fi}} 
              \def\@biblabel#1{#1.}
\makeatother

\date{March 2000}

\begin{document}

\title{Minimum $L_1$-Distance Projection onto\\
       the Boundary of a Convex Set: \\
       Simple Characterization%
       \footnote{This research was supported by CIBC Project No.~192251,
               funded jointly by the
               Canadian Imperial Bank of Commerce, and
               the Natural Sciences and Engineering Research Council of Canada.
               The author is grateful to Professor Wade~D. Cook and Dr.~Moshe Kress 
               for introducing the problem. \    }}
\author{Hans J.H. Tuenter%
        \footnote{Postdoctoral Fellow, Schulich School of Business, York University,
                  Toronto, Ontario, Canada.}}

\maketitle

\paragraph*{Abstract.}
We show that the minimum distance projection 
in the $L_1$-norm from an interior point onto the boundary of a convex set
is achieved by a single, unidimensional projection.
Application of this characterization when the convex set is a polyhedron
leads to either an elementary minmax problem or a set of easily solved linear programs,
depending upon whether the polyhedron is given as the intersection of
a set of halfspaces, or as the convex hull of a set of extreme points.
The outcome is an easier and more straightforward derivation of the ``special case''
results given in a recent paper by Briec~(Ref.~\cite{Briec:1997}).

\paragraph*{Key~Words.}
   Convex sets, minimum distance projection, $L_1$-norm.


\section{Introduction}
In a recent paper, Briec~(Ref.~\cite{Briec:1997}) studied the problem of 
determining the minimum distance from an interior point to the boundary of a convex set.
The approach taken by Briec is largely geometric and is based on the idea of separating convex sets with hyperplanes. This geometric viewpoint leads to a duality result, that is then applied to polyhedral sets.
For the special case that the distance is measured in the $L_1$-norm, Briec showed that 
the global minimum is obtained by solving a finite set of linear programs.

The purpose of this note is to show that, for the $L_p$-norm with $p$ less than or equal to one, there is a much 
simpler geometric characterization of the minimum distance projection from an interior 
point onto the boundary of a convex set.
Application of this characterization to polyhedral sets 
immediately renders Briec's ``special case'' results.

\section{Characterization of the Minimum $L_1$-Distance Projection}
Despite the vast literature on convex sets (cf.~Refs.~\cite{GruberWills:1993} and~\cite{Rockafellar:1970}), 
the following simple result involving the $L_1$-distance norm seems to have passed unnoticed.
\begin{theorem} \label{theo:1} \rm
 Let $S$ be a convex set in $\BBR^n$
 with the origin as an interior point,
 and let $X$ be the set of points on its boundary $\partial{S}$
 with minimum $L_1$-distance to the origin.
 Then there is an $x\in X$ with exactly one nonzero
 component.
\end{theorem}
\begin{proof}
 Let $x\in X$ be a point on the boundary of $S$ with minimum
 distance $d=\sum_i\left|x_i\right|$ to the origin.
 Note that $x$ has at least one nonzero component
 since the origin is an interior point.
 \par
 Now consider the set $E$, consisting of the $2n$ points $\pm de_i$,
 where $e_i$ is a vector with all components equal to zero
 except the $i$th component, which is equal to one.
 If there is an $\tilde{x}\in E$ that is not in the closure of $S$,
 then there is a $\lambda\in\left]0,1\right[$
 such that $\lambda\tilde{x}\in\partial{S}$.
 But then $d(\lambda\tilde{x})=\lambda d(\tilde{x})<d$,
 leading to a contradiction,
 and so $E\subset{\rm cl\,}S$ must hold.
 Now assume that all the elements of $E$ are interior points of~$S$.
 Any element of~$X$ can be expressed as a convex combination
 of elements of~$E$ as follows:
 $x=\sum_i (\left|x_i\right|/d)\{{\rm sgn}(x_i)de_i\}$,
 and this implies that $X\subset{\rm int\,}S$.
 This, now leads to a contradiction, and shows that
 there must be an $x\in E$ that lies on the boundary $\partial{S}$,
 and so completes the proof.
\end{proof}
This theorem implies that the minimum $L_1$-distance projection from an interior point 
to the boundary of a convex set is achieved by a single, unidimensional projection.
As the $L_1$-norm is invariant under translation,
the following corollary is immediate.
\begin{corollary} \label{corol:1}  \rm
 Let $a$ be an interior point of the convex set~$S$.
 If $d$ is the minimum $L_1$-distance from point~$a$ 
 to the boundary~$\partial{S}$, 
 then there is an $x\in\partial{S}$ with
 $d=\sum_i\left|x_i-a_i\right|$,
 that differs in only one component from~$a$.
\end{corollary}
%
Theorem~\ref{theo:1} and its corollary also hold when we take the $L_p$-norm with $p$ less than one 
as the distance metric. 
The proof proceeds along the same lines as for the case when $p$ is equal to one. 
The argument that any $x\in X$ can be expressed
as a convex combination of the elements of $E$ is a consequence of the fact that,
for $p\le1$,
\[ \max\left\{ \sum  \left|x_i\right|/d 
          \ \left| \, \vphantom{\sum} \right.
    \left(\sum\left|x_i\right|^p \right)^{1/p}=d 
\right\}
   = 1.
\]
For $p>1$, this maximum has a value of $n^{1-1/p}$,
and shows that [the proof of] Theorem~\ref{theo:1} does not generalize for the $L_p$-norm, 
when $p$ is larger than one. 

\section{Applications to Polyhedral Sets}
Briec considers two representational forms of a polyhedron:
first as the intersection of some finite collection of halfspaces,
and second as the convex hull of a finite set of extreme points.
For the former he remarks that, in the particular case of the $L_1$-metric, 
it is possible to find a global solution by solving $m$ linear programs,
where $m$ is the number of hyperplanes defining the polyhedron. 
In fact, as shown below, 
the global minimum can be obtained by solving an elementary minmax problem.
For the latter representation as a convex hull he shows that the global minimum 
is obtained as the minimum of $2m$ linear programs, where $m$ is the number of extreme points
defining the convex hull. 
In what follows, we assume for both representations that the origin is an interior point of the polyhedron.

\subsection{Intersection of Halfspaces}
If the polyhedron is given as the intersection of halfspaces it can be represented as
\[ S=\{x\in\BBR^n\mid Ax\le b\}.  \]
Note that the components of the vector~$b$ are strictly positive, as the origin is an interior point.
Because of this, we may assume, without loss of generality, that the matrix~$A$ does not contain a row that consists entirely of zeros, as this defines a redundant constraint.
By Theorem~\ref{theo:1} we know that there is a point on the boundary of $S$ with minimum distance to the origin  
that is a unidimensional projection of the form
$\lambda e_j$, with $\lambda a_j\le b$ and $\left|\lambda\right|$ maximal,
for some~$1\le j\le n$. Taking into account the sign of $\lambda$, define
\begin{eqnarray*}
  \lambda^-_j =  \max\, \left\{\lambda \mid -\lambda a_j\le b\right\}  
  \ \ \ {\rm and\ \ \ }
 \lambda^+_j =  \max\, \left\{\lambda \mid \lambda a_j\le b\right\}.  
\end{eqnarray*}
This gives the minimum distance as $d=\min_j\thinspace \{\lambda_j^-,\lambda^+_j\}$.
Writing out the constraints by component gives 
\[ \lambda_j^- = \min_{\{i\mid a_{ij}<0\}} -b_{i}/a_{ij} 
   \ \ \ {\rm and\ \ \ }
   \lambda_j^+ = \min_{\{i\mid a_{ij}>0\}} b_{i}/a_{ij}.
\]
We note that, in the case where an index set 
is empty,
which can occur when the polyhedron is unbounded,
one assigns to the corresponding $\lambda_j^-$ or $\lambda_j^+$ a value of infinity
or some sufficiently large number.
Now combine these expressions to 
\[ d=\min_i\, b_i/\max_j |a_{ij}|, \]
and one obtains the minimum distance as the solution to a simple minmax problem.
\subsection{Convex Hull}
If the polyhedron is given as the convex hull of a finite set of points
it can be represented as 
\[ S=\left\{x\in\BBR^n\mid x=\sum_k\nolimits  a_k \mu_k, \sum_k\nolimits\mu_k\le1, \mu_k\ge0\right\}. \]
Applying Theorem~\ref{theo:1}, while taking into account the sign of $\lambda_j$, gives
the following set of linear programs:
\begin{eqnarray*}
  \lambda^-_j &=& 
  \max\, \left\{\lambda \mid 
        -\lambda e_j= A\mu,\,
        \sum_k\nolimits\mu_k\le1,\,
        \mu_k\ge 0 
        \right\},  \\[0.5\baselineskip]
  \lambda^+_j &=& 
  \max\, \left\{\lambda \mid 
        \phantom{-} \lambda e_j=A\mu,\, 
        \sum_k\nolimits\mu_k\le1,\,
        \mu_k\ge0 
        \right\},  
\end{eqnarray*}
and the minimum distance as $d=\min_j\thinspace \{\lambda_j^-,\lambda^+_j\}$.
These linear programs are the duals of the ones in Briec's paper.

\section{Discussion}
Of course, our characterization and results do not have the generality that Briec's duality result (Proposition~3.1) has.
His result also covers the $L_p$-norm for $p>1$, in particular the important case of the Euclidean distance ($p=2$).
However, for the special cases considered here ($p\le1$), our characterization does allow for an easier geometric interpretation, and this leads to the corresponding optimization problems in a straightforward manner. 
Note however, that the characterization does not make any assumptions on the form of the convex set,
and therefore that it is applicable to convex sets other than polyhedral.

Finally, we note that a weighted version of the minimum distance projection,
where one considers the sum $\sum_i w_i\left|x_i\right|^p$, with the weights $w_i$ nonnegative 
and $p\le1$, 
also leads to the characterization that the minimum distance projection 
onto the boundary of the convex set
is achieved by a single, unidimensional projection.

\end{document}